\documentclass[a4paper,10pt]{article}
\usepackage{graphicx} 
\usepackage{amsmath} 
\usepackage{amssymb}  
\usepackage{amsfonts}
\usepackage{color}
\usepackage[english]{babel}
\usepackage[latin1]{inputenc}

\newcommand{\xmin}{x^{\star}}
\newcommand{\xtrue}{\widetilde{x}}

\newcommand{\dy}{\delta_y}
\newcommand{\dA}{\delta_A}
\newcommand{\DA}{\Delta_A}

\newcommand{\Dy}{\Delta_y}
\newcommand{\qA}{\mathcal{Q}(A)}
\newcommand{\qy}{\mathcal{Q}(y)}

\newcommand{\FS}{\mathcal{D}}

\newcommand{\R}{\mathbb{R}}

\newcommand{\one}{\mathbf{1}}

\newtheorem{proposition}{Proposition}

\newtheorem{assumption}{Assumption}
\newenvironment{proof}{\paragraph{Proof:}}{\hfill$\square$}

\usepackage{graphics} 
\usepackage{epsfig} 

\usepackage{mathtools}

\usepackage{tikz}
\usetikzlibrary{shapes,arrows,calc,positioning}

\tikzset{
	block/.style    = {draw, thick, rectangle, minimum height = 2.5em,
		minimum width = 2.5em, node distance = 1.7cm},
     sum/.style      = {draw, circle, radius = 0.5cm, node distance = 1.3cm}, 
     input/.style    = {coordinate}, 
     output/.style   = {coordinate} 
 }

{\normalfont}{\normalfont}
\newtheorem{result}{\noindent \textbf{Result}}{\normalfont}{\normalfont}
 {\normalfont }{\normalfont}
{\normalfont}{\normalfont}

\begin{document}
\title{Sparse linear regression with compressed and low-precision data via concave quadratic programming}
\author{V. Cerone, S. M. Fosson$^{*}$, D. Regruto 
}

\maketitle
\begin{abstract}
We consider the problem of the recovery of a $k$-sparse vector from compressed linear measurements when data are corrupted by a quantization noise. When the number of measurements is not sufficiently large, different $k$-sparse solutions may be present in the feasible set, and the classical $\ell_1$ approach may be unsuccessful. For this motivation, we propose a non-convex quadratic programming method, which exploits prior information on the magnitude of the non-zero parameters. This results in a more efficient support recovery. We provide  sufficient conditions for successful recovery and numerical simulations to illustrate the practical feasibility of the proposed method.
\end{abstract}

\section{Introduction}\label{sec:intro}

Sparse linear regression is the recovery of a sparse vector $x\in\R^n$ from linear measurements $y=Ax\in\R^m$, with $A\in\R^{m,n}$. A vector is sparse if it has few non-zero components; more precisely, we call $k$-sparse a vector with $k\ll n$ non-zero components. The interest for sparse solutions has different motivations. In machine learning and system identification, a purpose is to build models as simple as possible from large datasets. Indeed, we know that in many cases the true number of parameters of a system is much smaller than the global dimensionality of the problem, and sparsity supports the removal of redundant parameters. In the recent literature, the identification of linear systems under sparsity constraints is considered in \cite{san11,tot11,tothsysid12,fox18cdc}. Furthermore, in the last decade, sparsity has been attracting a lot of attention due to the theory of compressed sensing (CS, \cite{don06,fou13}), which states that a sparse vector can be recovered from compressed linear measurements, that is, when $m<n$, under suitable conditions on the structure of $A$.

In real applications, observations may be affected by different sources of noise. The most common one is measurement noise, which is usually modeled as an  unknown additive term $\dy\in\R^m$, i.e., $y=Ax+\dy$.  Linear regression in the presence of measurement noise is usually solved via least squares; if $x$ is desired to be sparse, an $\ell_1$ regularizer is added, which promotes sparsity and keeps the problem convex. More precisely, the problem can be formulated as $\min_{x\in\R^n}\|x\|_1$ s. t. $\|Ax-y\|_p\leq \epsilon$, where $\|\cdot\|_p$ is a suitable norm and $\epsilon>0$ is a known bound; this is referred to as Basis Pursuit Denoising (BPDN$_p$). In the literature,  BPDN$_2$ is very popular, in particular in signal processing and CS, see, e.g., \cite{fou13}. The choice $p=2$ bounds the mean energy of the error, hence BPDN$_2$ is tolerant to possible outliers.  The case $p=\infty$, instead, is considered to deal with quantized or low-precision data in \cite{don06infty,las11,val15,fox19acc}. When $y$ is quantized, a bound on each component is given, and no outliers occur. This makes  the $\ell_{\infty}$ description more suitable than the $\ell_2$ one.

In many applications, both the measurements and the matrix $A$ are affected by noise. In particular, in linear systems identification, this is analyzed by the errors-in-variables (EIV) model, where both the input and the output are corrupted by noise; we refer the reader to \cite{sod18} for a complete overview on EIV models. In \cite{cer93,cer12}, the identification of the feasible set for EIV models is tackled. Perturbations on $A$ are more tricky to deal with than measurement noise, as they make the problem non-convex. Furthermore, in CS, the recovery of sparse vectors from compressed measurements is strongly compromised by perturbations on $A$, which are tolerated only if sufficiently small or with particular structures. In \cite{her10}, robustness to perturbation in CS is investigated, providing information on the amount of perturbation on $A$ that can be tolerated by BPDN$_2$.  In \cite{yan12}, the perturbation is assumed to be $n$-dimensional, instead of $nm$-dimensional, and bounded in the $\ell_{\infty}$ norm. Alternating minimization is used to find a minimum, which is global when the components of $x$ are known to be non-negative, see \cite[Section IV.A]{yan12}.

In \cite{fox19acc}, the same non-negativity assumption is leveraged to formulate the sparse linear regression problem, with bounded noise both on $y$ and $A$, as a  linear programming problem. Specifically, the proposed problem in \cite{fox19acc} is the minimization of the $\ell_1$ norm within the feasible set, which, under the non-negativity assumption, is a convex polytope. The approach is proved to be robust, in the sense that the distance between the solution and the true vector is proportional to the noise and to the sparsity level of the true signal. 

The $\ell_1$ approach may fail if the feasible set does not delineate the sparsity pattern. More precisely, many $k$-sparse solutions might belong to the feasible set, though with different support.  When several $k$-sparse solutions are feasible, the $\ell_1$ minimizer selects the one with less energy, which may not correspond to the desired one. We highlight that, particularly in the case of compressed measurements, the support recovery is the key point for the vector recovery. If the support is known, and if we assume sparsity level $k\leq m$, the problem is not underdetermined, and the solution can be found by inversion, see, e.g., \cite{fox14}. We also remark that, in CS, conditions on $A$ are provided that guarantee a unique $k$-sparse solution for BPDN$_2$, see, e.g., \cite{fou13}; nevertheless, these conditions are actually limited to specific classes of random matrices, and perturbations on $A$ are not envisaged.

In this work, we tackle the above described problem by proposing a novel solution. Specifically, we propose a non-convex, polynomial approach which is oriented to improve the recovery performance,  particularly in terms of sparsity pattern, by leveraging a prior information on the magnitude of the non-zero components. Sufficient conditions that guarantee the recovery of the correct support are proven, and numerical simulations illustrate the effectiveness of the approach in different settings.

The paper is organized as follows. In Section \ref{sec:me}, we present the background  and a motivating example, which highlights the flaws of the $\ell_1$  approach. In Section \ref{sec:ps}, we formally state the proposed problem, and discuss the considered assumptions. In Section \ref{sec:pop}, we explain how to solve the proposed problem. In Section \ref{sec:an}, we prove sufficient conditions that guarantee successful recovery, while Section \ref{sec:nr} is devoted to numerical experiments. Finally, some conclusions are drawn in Section \ref{sec:con}.
\section{Background and motivating example}\label{sec:me}
Let us consider an over-dimensioned vector of parameters $\xtrue\in\R^n$, i.e., only $k\ll n$ components are different from zero, and let us assume that the system output is given by $y=A\xtrue$, $A\in\R^{m,n}$. The goal is to recover $\xtrue$ under the following assumptions: (a) $m<n$, i.e., compressed output observations are available; (b) both $y$ and $A$ are not exactly known. More precisely, low-precision (or quantized) versions of $A$ and $y$ are known, which are denoted by $\qA\in\R^{m,n}$ and $\qy\in\R^m$, respectively. The corresponding perturbations are denoted by $\dA=\qA-A$ and $\dy=\qy-y$. Moreover, we assume to know the maximum perturbation magnitudes, denoted by $\DA$ and $\Dy$, respectively. For example, if we consider data quantized over a fixed number of bits, we assume to know the number of bits.

Let us represent by $\FS$ the feasible set. The natural approach to estimate $\xtrue$, which is sparse, is by looking for the sparsest solution in $\FS$. This can be performed by minimizing the $\ell_1$-norm over $\FS\subset \R^n$, as illustrated in \cite{fox19acc}. Specifically, the problem is formulated as follows: given $\qy\in\R^m$, $\qA\in\R^{m,n}$, $\DA>0$, and $\Dy>0$,
\begin{equation}\label{theproblem1}
 \begin{split}
  \min_{x\in\R^n}\|x\|_1~\text{ s. t. } &y=A x\\
  &\qy=y+\dy\\
  &\qA=A+\dA\\
  &\|\dy\|_{\infty}\leq \Dy\\
  &\|\dA\|_{\infty}\leq \DA.
 \end{split}
\end{equation}
In \cite{fox19acc}, the assumption $x_i\geq 0$, for any $i=1,\dots,n$, is considered, which makes the problem solvable via linear programming. Actually, the same result is obtained if each component is known to be either non-positive or non-negative: non-negativity is assumed for simplicity without loss of generality. However, we notice that the non-negative setting has itself a number of applications, that range from sparse localization problems, see \cite{bay15}, to sparse sensor selection, see \cite{cal16}.

In \cite{fox19acc}, Problem \eqref{theproblem1} is shown to provide encouraging results if compared to state-of-the-art algorithms for the recovery of sparse signals from compressed and low-precision data.
\begin{figure}
\centering
\includegraphics[width=0.4\textwidth]{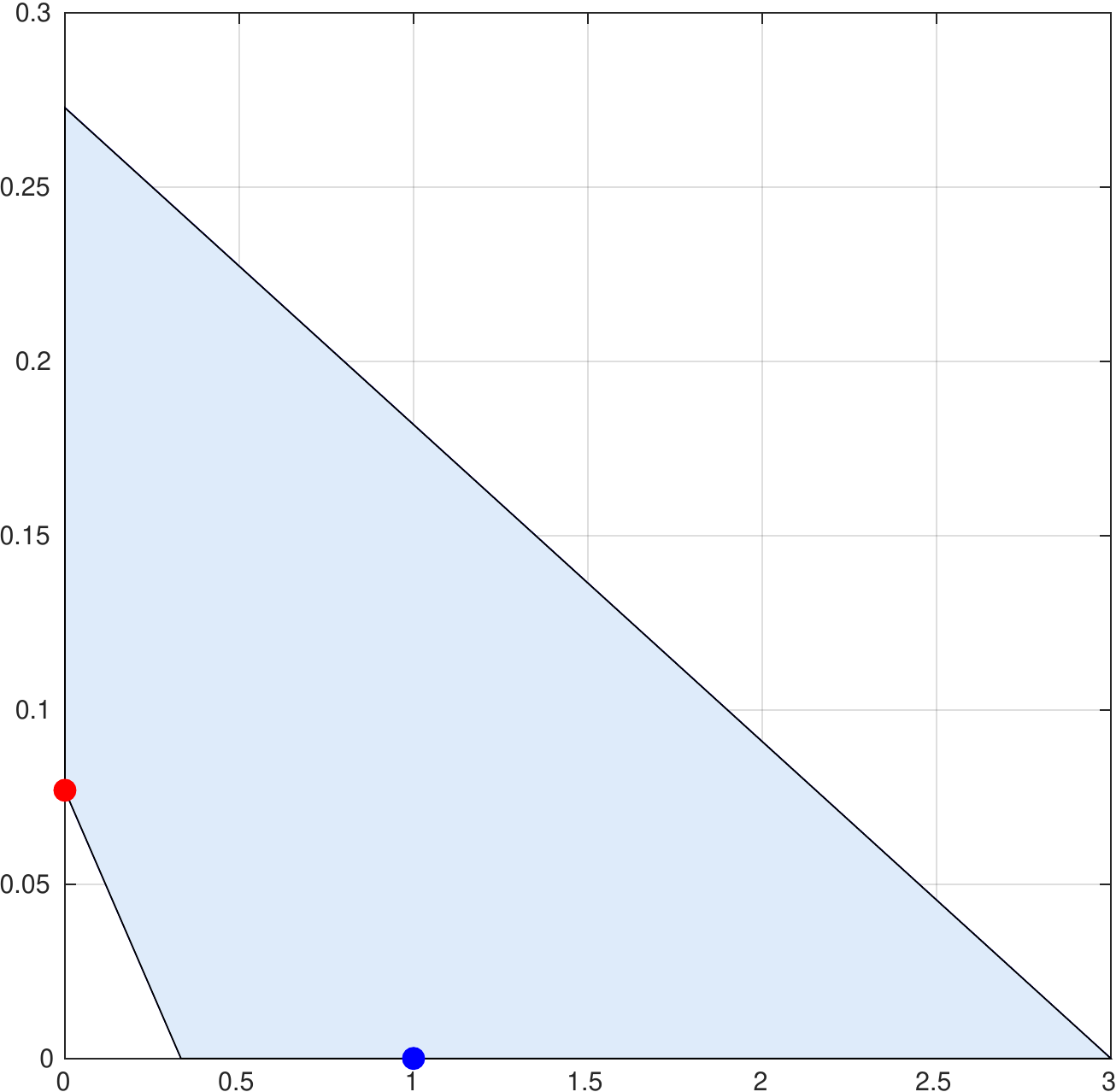}
\caption{$A=[0.2131,1.2414]$, $\xtrue=[1,0]^T$ (blue point); $\qA=[0.2,1.2]$, $\qy=0.2$. The feasible set $\FS$ is the light blue polytope. The solution of Problem \eqref{theproblem1}  is $[0,0.0769]^T$ (red point). The support is inverted with respect to $\xtrue$.}\label{fig0}
\end{figure}
However, a drawback of the $\ell_1$ approach is that it may fail the recovery of the right solution when several $k$-sparse vectors are present in $\FS$, with different supports. We illustrate this occurrence with a simple motivating example.
Let us consider $A=[0.2131,1.2414]$, $\xtrue=[1,0]^T$. Then we have $y=A\xtrue=0.2131$. Data are subject to quantization with step $0.1$, thus we observe $\qA=[0.2,1.2]$ and $\qy=0.2$. If we assume non-negativity, the feasible set $\FS$, which is shown in Fig. \ref{fig0}, is the intersection between the first orthant and the area between the two lines that represent the constraints. The solution of problem \eqref{theproblem1} is $[0,0.0769]^T$, depicted in Fig. \ref{fig0} with a red point. This vector is as sparse as $\xtrue$, while it has the opposite support. Moreover, this solution has smaller non-zero value than the true one, which is not a natural choice in practice. In fact, in real applications, the null entries of $\xtrue$ are not exactly zero, instead they take small values; from this point of view, the support discriminates between significant values and non-significant ones. In our example, $0.0769$ is rather close to zero if compared to 1, therefore a suitable algorithm should estimate it as a zero, and search for a more significant non-zero value to put in the support.

To tackle this occurrence, intuitively, one could assume that an interval $[\alpha,\beta]$, $\alpha>0$, $\beta\geq \alpha$, of significant magnitudes is known, thus the constraint $x_i\in\{0\}\cup [\alpha,\beta]$, $i=1,\dots,n$, might be added to Problem \eqref{theproblem1}. Nevertheless, this would compromise the convexity of $\FS$  and make the problem definitely more complex. For this motivation, in this work we propose a different solution.

\section{Problem Statement}\label{sec:ps}
In this  section, we present a novel formulation of the sparse linear regression problem with compressed measurements and low-precision data, under the assumption that the significant values are in $[\alpha,\beta]$. We start by studying the case of non-negative parameters, while we leave for future work more general sign settings.

\begin{assumption}\label{ass2}
For any $i=1,\dots,n$, $x_i\in\{0\}\cup [\alpha,\beta]$.\\ 
\end{assumption}
Then, we define $d:=\frac{\alpha+\beta}{2}$, and we formulate the problem as follows.\\

\textbf{Problem 1}
\\
Given $\qy\in\R^m$, $\qA\in\R^{m,n}$, $\DA>0$, and $\Dy>0$,
\begin{equation}\label{theproblem}
 \begin{split}
  \min_{x\in [0,d]^n} d\|x\|_1-\|x\|_2^2&~\text{ s. t. } y=A x\\
  &\qy=y+\dy\\
  &\qA=A+\dA\\
  &\|\dy\|_{\infty}\leq \Dy\\
  &\|\dA\|_{\infty}\leq \DA.\\
 \end{split}
\end{equation}

Intuitively, Problem 1 is more accurate than \eqref{theproblem1} since the concave penalty $d\|x\|_1-\|x\|_2^2$ is closer to $\|x\|_0$ than $\|x\|_1$, see, e.g., \cite{fox18asi} for more details. The rest of the paper is devoted to prove this enhancement in the proposed setting.

First of all, we observe that if $x\in[0,d]^n$, then the minimum of the objective functional is zero, which is achieved for any $x\in\{0,d\}^n$. Since the true non-zero parameters are in $[\alpha,\beta]$, this approach does not exactly evaluate the non-zero parameters. However, our main goal is to obtain the right support. As mentioned before, once the right support is obtained, if $k\leq m$, the problem of evaluating the non-zero parameters is not underdetermined and can be obtained by inversion. Nevertheless, in this work, we assume $[\alpha,\beta]$ small enough such that  $d$ itself is a good approximation; the case of large intervals is left for future work.

In the next section, we illustrate how to practically tackle Problem 1 through polynomial optimization.

\section{Polynomial optimization approach}\label{sec:pop}
Thanks to the following result, we show that, under Assumption \ref{ass2}, Problem 1 can be solved through non-convex quadratic programming. 
In the following, given two vectors $a,b\in\R^n$, we write $a\succeq b$  to denote $a_i \geq b_i$ for each $i=1,\dots,n$. We denote by $I_n\in\R^{n,n}$ the identity matrix. Moreover, $\one_n:=(1,1,\dots,1)^T\in\R^n$. 
\\
\begin{result}
Under Assumption \ref{ass2}, Problem 1 can be equivalently written as the following polynomial problem:
\begin{equation}\label{sys3}
 \begin{split}
  \min_{x\in[0,d]^n}d&\|x\|_1-\|x\|_2^2~\text{ s. t. } C x \preceq g\\
  &\text{where }\\
  &C= \left(\begin{array}{c}
         \qA-\DA\one_m\one_n^T\\
         -\qA-\DA\one_m\one_n^T\\
        \end{array}\right)\in\R^{2m,n}\\
    &g= \left(\begin{array}{c}
         \qy+\Dy\one_m\\
         -\qy+\Dy\one_m\\
        \end{array}\right)\in\R^{2m}.\\
        \end{split}
\end{equation}
\end{result}

Result 1 is obtained by observing that Problem 1 can be rewritten as follows:
\begin{equation}\label{sys2}
 \begin{split}
  \min_{x\in[0,d]^n}d\|x\|_1-\|x\|_2^2\text{ s. t. } &\|\qy-(\qA-\dA)x\|_{\infty}\hspace{-0.05cm}\leq\hspace{-0.05cm}\Dy\\
  &\|\dA\|_{\infty}\leq \DA
 \end{split}
\end{equation}
and by applying the results on bounded EIV identification proved in \cite{cer93}, under Assumption \ref{ass2}.

We notice that Problem 1, as written in \eqref{sys3}, is non-convex and semi-algebraic, with linear constraints. In particular, when $x$ is known to have only non-negative entries, the $\ell_1$ term is a sum, and the problem is a concave quadratic programming problem. To solve Problem 1, we can apply polynomial optimization results by \cite{las01}, which state that the global minimum of a constrained polynomial problem can be achieved by a hierarchy  of relaxed semidefinite programming (SDP) problems. As the order of relaxation increases, the solutions of these SDP converge to the global optimal solution of the original problem. The hierarchy has finite convergence generically, as illustrated in \cite{nie13,nie14}.

The SDP relaxation technique is computationally intense in case of large dimensional data. However, in \cite{lasspa} it is shown that if the polynomial problem has a sparse structure, the dimension of the SDP relaxations can be reduced. This polynomial sparsity arises, for example, when the variables are decoupled, which is the case of the functional $d\|x\|_1-\|x\|_2^2$. Moreover, since the functional is quadratic, also the relaxation order can be kept low, as described in the numerical experiments in Section \ref{sec:nr}.

We remark that the generalization to a formulation that envisages unknown signs can be derived by exploiting results in \cite{cer93}. This generalization produces a feasible set which is not overall convex, though convex over each orthant.  The solution is then feasible, while computationally intensive for large dimensional problem. The analysis of this problem is left for future work.

\section{Analysis}\label{sec:an}
In this section, we provide conditions under which Problem 1 is successful when the non-zero parameters of the system are known to be in the range $[\alpha,\beta]$, $\alpha>0$, $\beta\geq \alpha$. We start from the simplified case where $\alpha=\beta=d$, i.e., $\xtrue\in\{0,d\}^n$.

The following proposition analyses conditions under which $\xtrue$ is the unique feasible solution in $\{0,d\}^n$.
\begin{proposition}\label{prop1}
Let $A\xtrue=y$, with $\xtrue\in\{0,d\}^n$. Let $A_i$ be the $i$-th column of $A$. If 
\begin{equation}\label{condition1}
\left\|\sum_{i=1}^n \gamma_i A_i\right\|_{\infty}>\frac{2\Dy}{d}, \text{ for any non-null } \gamma\in\{0,\pm 1\}^n,
\end{equation}
then $\xtrue$ is the unique solution of Problem 1.
\end{proposition}

\begin{proof}
Let us consider any $z\in\{0,d\}^n$, $z\neq \xtrue$. Then, $w:=\xtrue-z\in\{0,\pm d\}^n$. If $\|\sum_{i=1}^n \gamma_i A_i\|_{\infty}>\frac{2\Dy}{d}$ for any   non-null $\gamma_i\in\{0,\pm 1\}$, then 
\begin{equation}
\| A w\|_{\infty}> 2\Dy. 
\end{equation}
This implies:
\begin{equation}\label{calcoli}
\begin{split}
&\|A w \pm \dy \|_{\infty}> 2\Dy\\
\Rightarrow~&\|A w + \dy \|_{\infty}+\|\dy \|_{\infty}> 2\Dy\\
\Rightarrow~&\|A w + \dy \|_{\infty}> 2\Dy-\|\dy \|_{\infty}\geq \Dy.
\end{split}
\end{equation}
Since $w=\xtrue-z$ and $\qy=A\xtrue+\dy$, we obtain:
\begin{equation}\label{altricalcoli}
\begin{split}
&\|\qy-Az \|_{\infty}> \Dy.
\end{split}
\end{equation}
The last inequality follows from the fact that $\|\dy \|_{\infty}\leq \Dy$. 
This proves that $z$ does not belong to the feasible set. Therefore, $\xtrue$ is the unique feasible solution.
\end{proof}
Concerning the condition $\|\sum_{i=1}^n \gamma_i A_i\|_{\infty}>\frac{2\Dy}{d}$ for any  non-null $\gamma\in\{0,\pm 1\}^n$, we highlight that, as discussed in \cite{tib13}, in many applications the columns of $A$ are in general position,  i.e., $\sum_{i=1}^n \gamma_i A_i\neq 0$  for any  non-null  $\gamma\in\{0,\pm 1\}^n$. This definition is rather technical; however, it is proven to hold in many common cases, for instance, when $A$ has random entries generated from a continuous distribution, see Lemmas 3-4 in \cite{tib13}. Under the general position assumption, there exists a $\sigma>0$ such that $\|\sum_{i=1}^n \gamma_i A_i\|_{\infty}\geq \sigma$. Therefore, if $d$ is sufficiently large with respect to $\Dy$, we can expect $\frac{2\Dy}{d}<\sigma$. On the other hand, $\sigma$ can be increased by increasing $m$.

Now, we extend Proposition \ref{prop1} to the case of non-null $\xtrue_i$ belonging to $[\alpha,\beta]$.

\begin{proposition}\label{prop2}
Let $A\xtrue=y$, with $\xtrue_i \in \{0\}\cup [\alpha,\beta]$ for any $i=1,\dots,n$. Let $A_i$ be the $i$-th column of $A$. If 
\begin{equation}\label{condition2}
\left\|\sum_{i=1}^n \gamma_i A_i\right\|_{\infty}>2\Dy \text{ for any  non-null } \gamma\in Q^n,
\end{equation}
where 
$$Q:=\{-d\}\cup \left[\frac{\alpha-\beta}{2},\frac{\beta-\alpha}{2}\right]\cup \left[\alpha,\beta\right]$$ and $d=\frac{\alpha+\beta}{2}$, then Problem 1 has a unique solution $\xmin\in\{0,d\}^n$, and $\xmin$ has the same support of $\xtrue$.
\end{proposition}

\begin{proof}
We denote by $S(v)$ the support of any $v\in\R^n$. 
Let us consider any $z\in\{0,d\}^n$, $S(z)\neq S(\xtrue)$. Let $w:=\xtrue-z$. Then, $w\in Q^n$.

If $\|\sum_{i=1}^n \gamma_i A_i\|_{\infty}> 2\Dy$ for any   $\gamma\in Q^n$, then  
\begin{equation}
\| A w\|_{\infty}> 2\Dy. 
\end{equation}

By computations similar to \eqref{calcoli} and \eqref{altricalcoli}, we can conclude that $z$ is not in the feasible set. Therefore, the solution of Problem 1 is the unique $\xmin\in\{0,d\}^n\cap \FS$ such that $S(\xmin)=S(\xtrue)$.
\end{proof}

We finally notice that in Proposition \ref{prop1} and Proposition \ref{prop2}, the condition  $\|\sum_{i=1}^n \gamma_i A_i\|_{\infty}>\sigma$ for some $\sigma>0$ is defined on $A$, which is not known. However, even if unknown, it makes sense to assume that some information on the properties of $A$ is available. In particular, the fact that the general position of the columns holds can be assumed in many applications. 
Nevertheless, in the case that $\sigma$ has to be assessed more precisely, one can state the condition on the observed $\qA$, and the following sufficient condition can be proven.
\begin{proposition}\label{prop3}
Let $A\xtrue=y$, with $\xtrue_i \in \{0\}\cup [\alpha,\beta]$ for any $i=1,\dots,n$. Let $\qA=A+\dA$, and let $\qA_i$ be the $i$-th column of $\qA$. If 
\begin{equation}\label{condition3}
\left\|\sum_{i=1}^n \gamma_i \qA_i\right\|_{\infty}>2\Dy+\DA\beta n~ \text{ for any non-null } \gamma\in Q^n
\end{equation}
where $Q=\{-d\}\cup \left[\frac{\alpha-\beta}{2},\frac{\beta-\alpha}{2}\right]\cup \left[\alpha,\beta\right]$, $d=\frac{\alpha+\beta}{2}$, then, Problem 1 has a unique solution $\xmin\in\{0,d\}^n$, and $\xmin$ has the same support of $\xtrue$.
\end{proposition}
The proof is similar to that of Proposition \ref{prop2}, thus omitted for brevity. As expected, the sufficient condition in Proposition \ref{prop3} is more restrictive than that of Proposition \ref{prop2}, while the two tend to coincide when the perturbation on $A$ tends to zero.

If we consider the motivating example illustrated in Section \ref{sec:me}, condition \eqref{condition1} is fulfilled if $d=1$, since $|A\gamma|\geq 0.2131>0.2=2 \Dy$. Therefore, our approach provides the right solution. Furthermore, condition \eqref{condition2} is fulfilled, for example, with $[\alpha,\beta]=[0.9,1.1]$ and $\Dy=0.03$.
\section{Numerical results}\label{sec:nr}
In this  section, we present numerical results, which illustrate the effectiveness of the proposed approach with respect to the $\ell_1$ approach in \cite{fox19acc}. We show two experiments. In the first one, we consider the ideal setting where the true vector of parameters is $\xtrue\in\{0,d\}^n$; thus, in this case, we expect to recover exactly both the support and the non-zero values, under suitable conditions. We specify that, despite its peculiarity, the binary setting is relevant in widespread applications, such as localization and image processing, see \cite{fox18asi,fox19abu} for a complete overview. In the second experiment, instead, the non-zero components are in $[\alpha,\beta]$, and $d=\frac{\alpha+\beta}{2}$ ; in this case, we expect to obtain the right support and a biased estimation of the non-zero values.

\begin{figure*}
\centering
\includegraphics[width=0.46\textwidth]{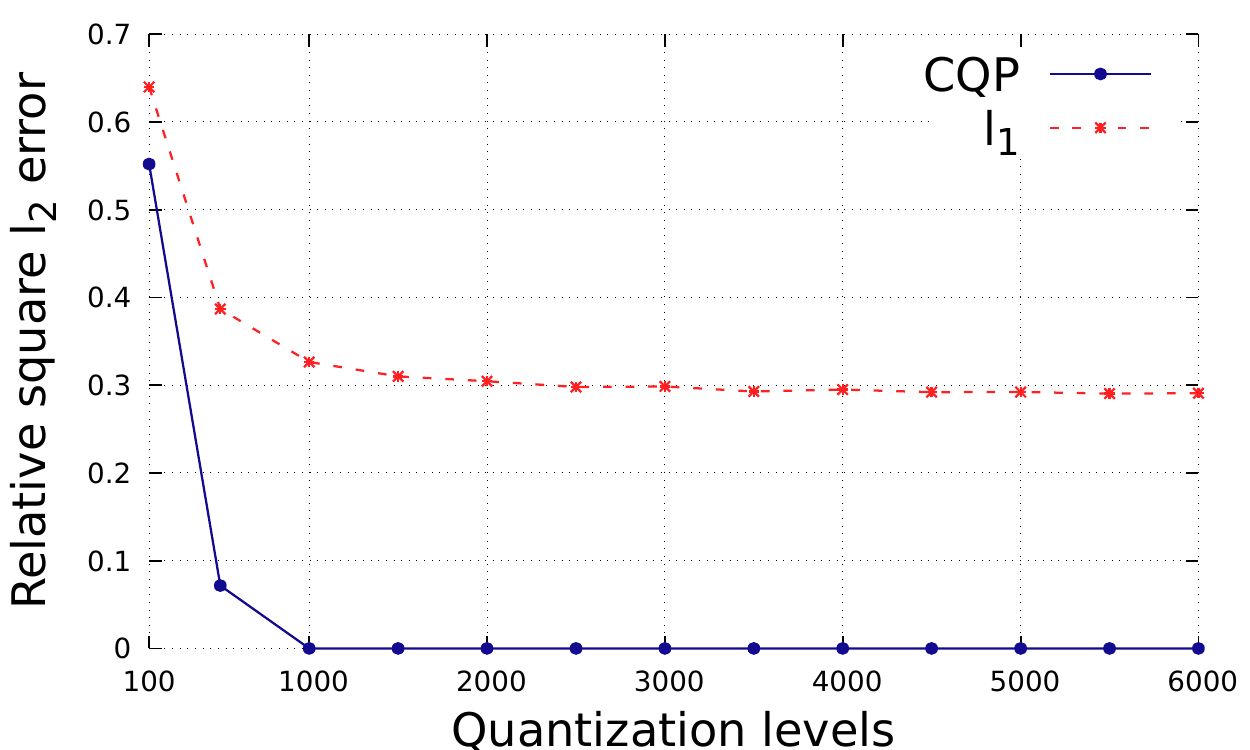}\quad
\includegraphics[width=0.46\textwidth]{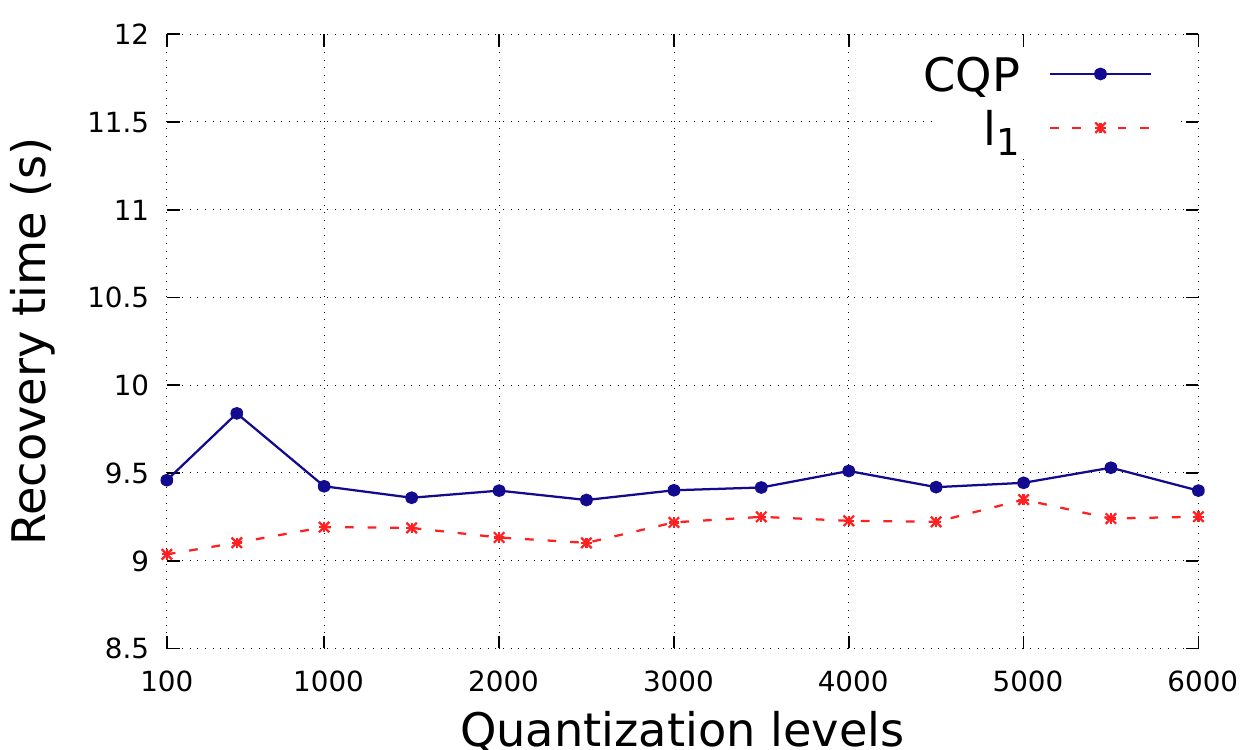}
\includegraphics[width=0.46\textwidth]{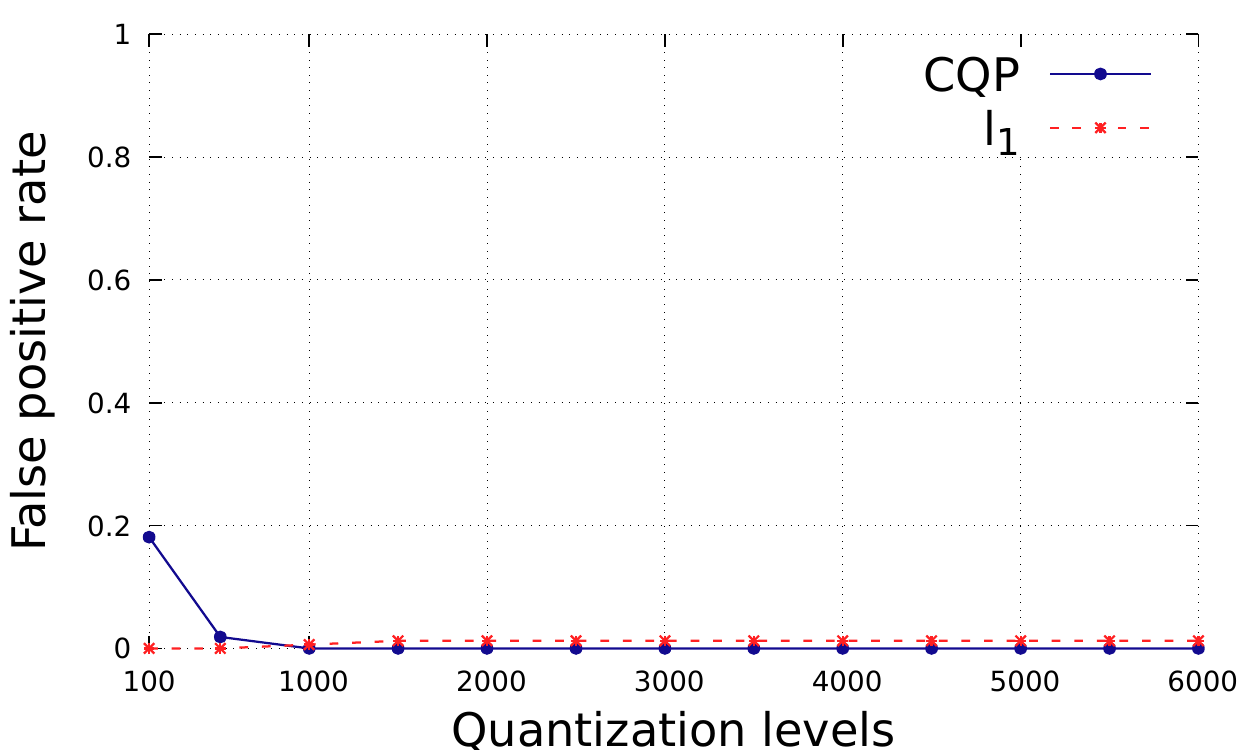} \quad
\includegraphics[width=0.46\textwidth]{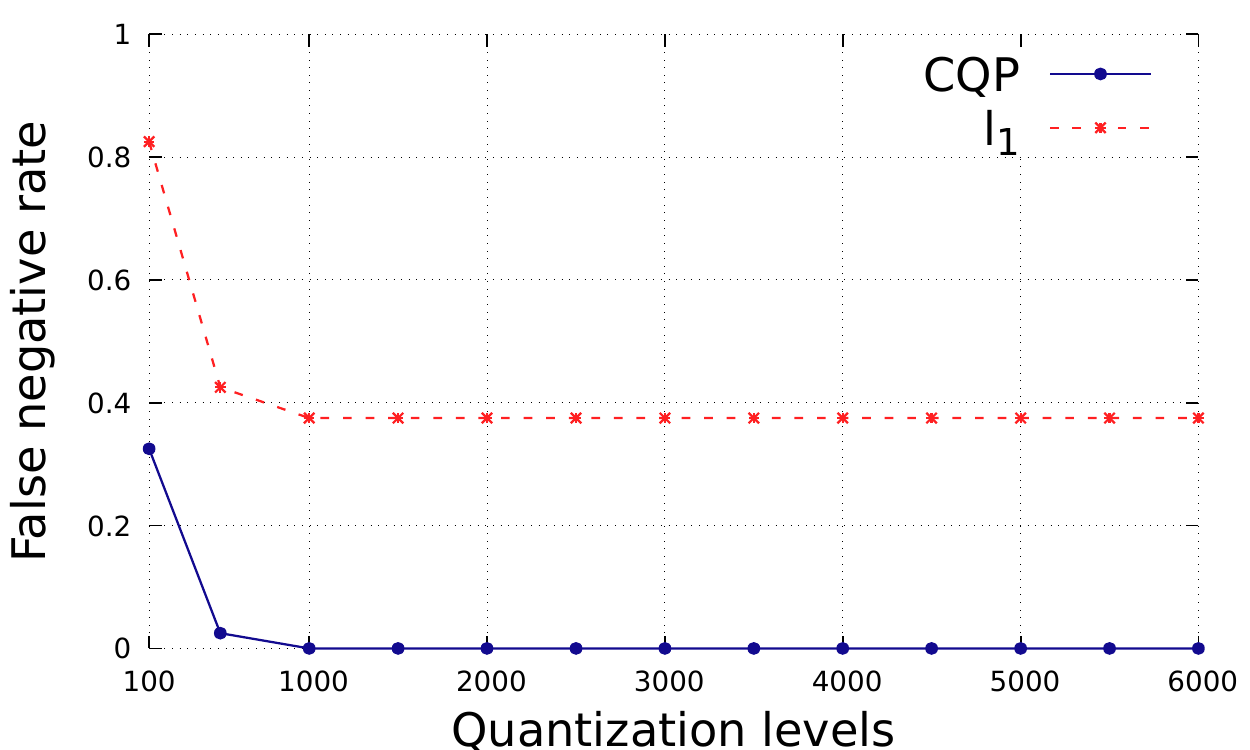}
\caption{Experiment 1: comparison between the proposed method, denoted by CQP, and the $\ell_1$ method \cite{fox19acc}, under the assumption that $\xtrue\in\{0,d\}^n$.}\label{fig1:errors}
\end{figure*}
\begin{figure*}
\centering
\includegraphics[width=0.46\textwidth]{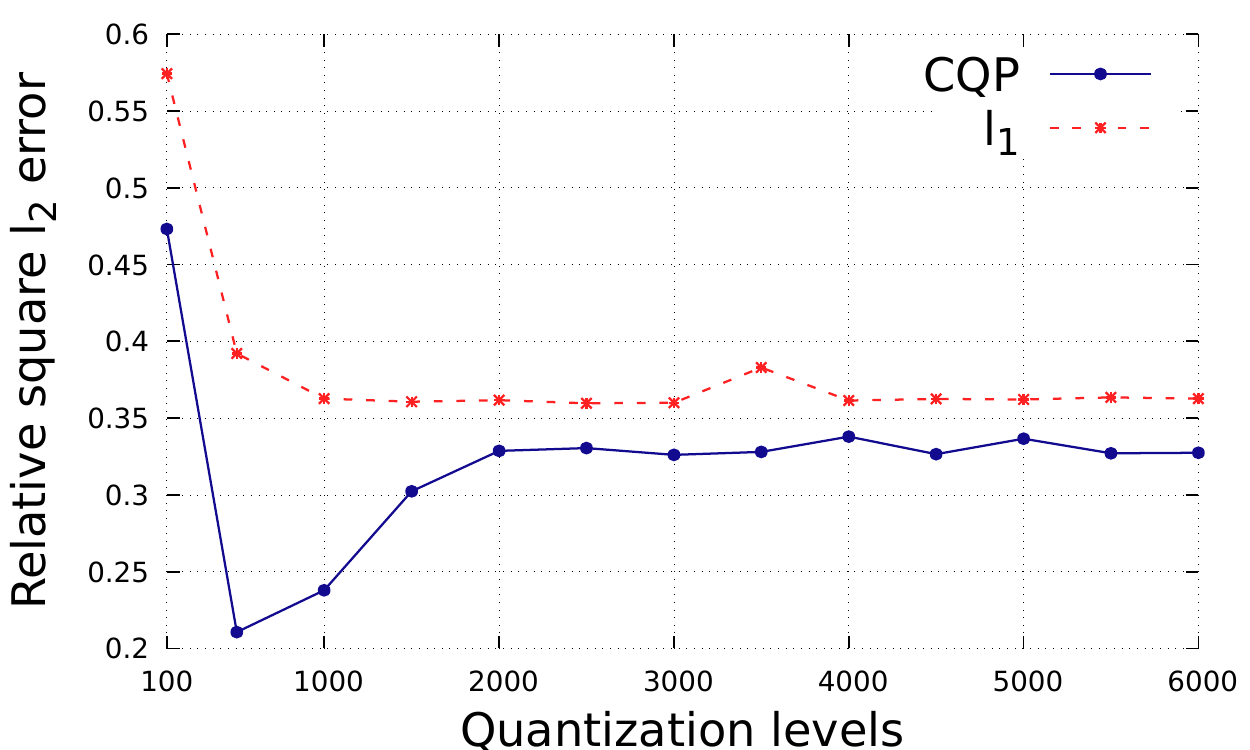}\quad
\includegraphics[width=0.46\textwidth]{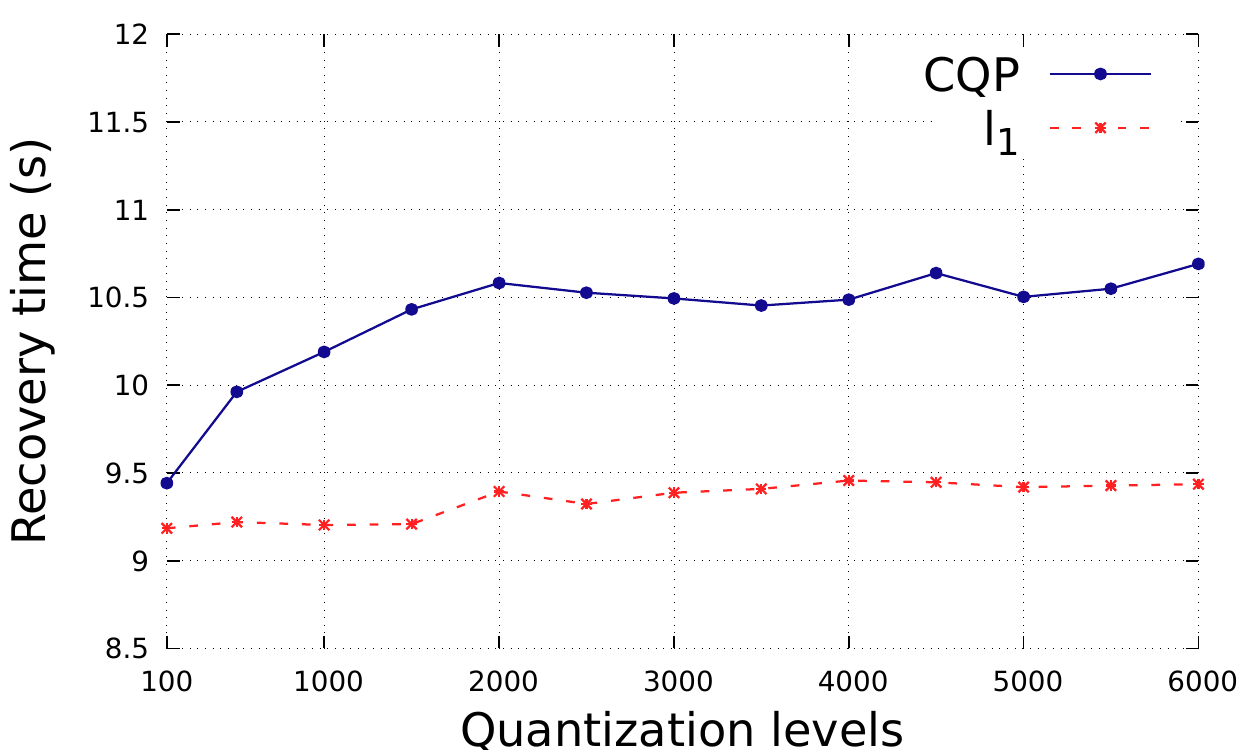}
\includegraphics[width=0.46\textwidth]{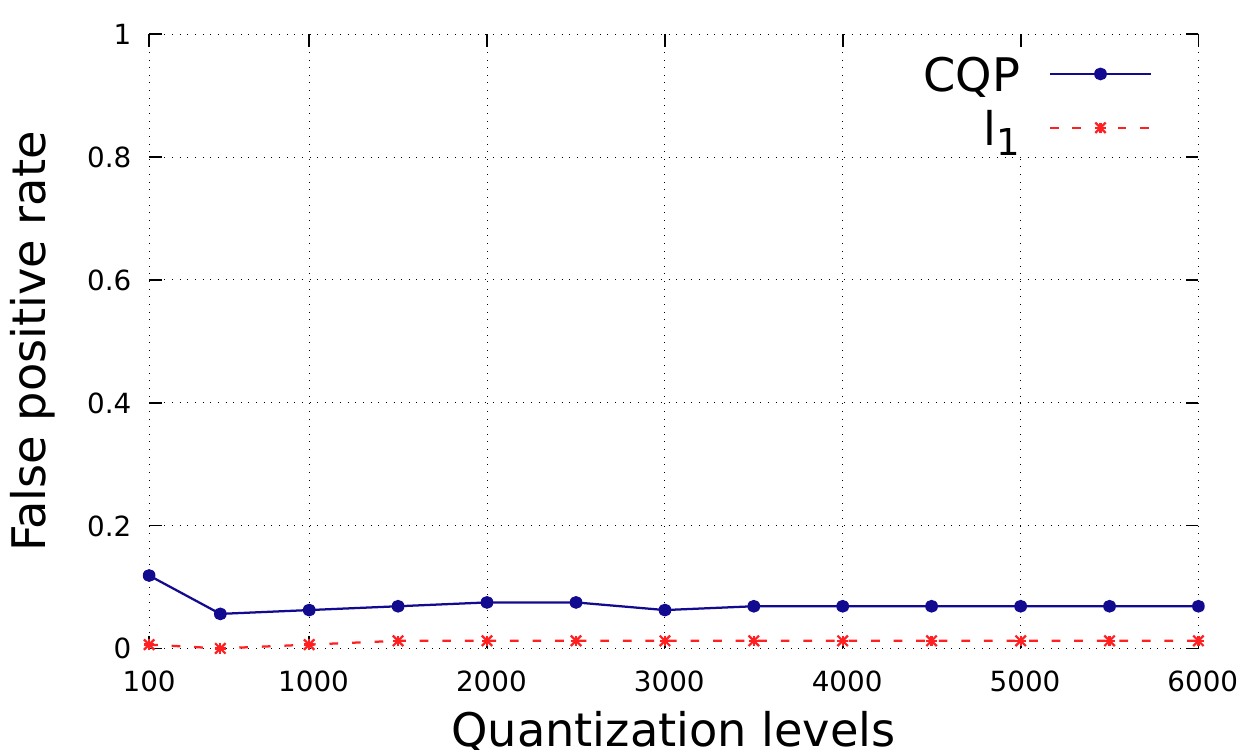} \quad
\includegraphics[width=0.46\textwidth]{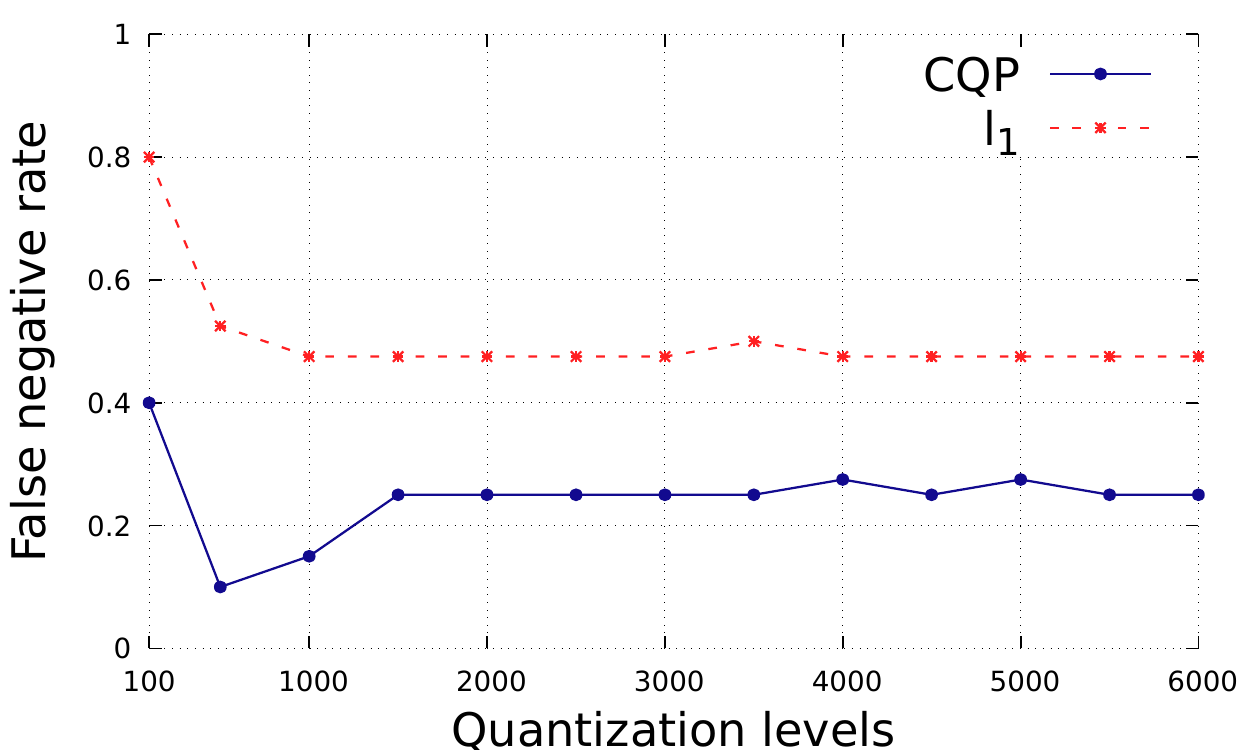}
\caption{Experiment 2: comparison between the proposed method, denoted by CQP, and the $\ell_1$ method \cite{fox19acc}, under the assumption that the non-zero values of $\xtrue$ are in $[0.8,1.2]$.}\label{fig2:errors}
\end{figure*}
For both experiments, we propose the following setting. A system acquires a sparse vector $\xtrue\in\R^n$ , with $n=10$, through a regression matrix $A\in\R^{m,n}$, $m=4$. The sparsity level of $\xtrue$ is $k=2$, i.e., only two parameters are non-null. The knowledge of $k$ is not required by the considered  recovery algorithms. $A$ is generated according to a Gaussian distribution $\mathcal{N}(0,\frac{1}{m})$; the support of $\xtrue$ is generated uniformly at random, while the non-zero entries are uniformly distributed in $[\alpha,\beta]$. We assume that data are uniformly quantized: given a certain number of equidistant quantization levels, each entry of $A$ and $y$ is approximated with the closest point in the quantization codebook. The quantization range is considered  sufficiently large to avoid saturation problems. The considered quantization levels range from 100 to 6000, both for $A$ and $y$.

In the first experiment, depicted in Fig. \ref{fig1:errors} we consider the case $d=\alpha=\beta=1$; in the second experiment, illustrated  in Fig. \ref{fig2:errors}, $[\alpha,\beta]=[0.8,1.2]$. The parameters are supposed to be non-negative. The proposed method, based on concave quadratic programming (CQP in the figures) is compared to the $\ell_1$ approach of \cite{fox19acc} in terms of different performance metrics: 
the relative square $\ell_2$ error, defined as $\|\widehat{x}-\xtrue\|_2^2/\|\xtrue\|_2^2$, where $\widehat{x}\in\R^n$ is the estimate; the false positive rate, defined as the number of events where $\widehat{x}_i\neq 0$ while $\xtrue_i= 0$, over $n-k$; the false negative rate,  defined as the number of events where $\widehat{x}_i = 0$ while $\xtrue_i\neq 0$, over $k$; the run time (in seconds).
Experiments are performed on Matlab R2016b, on a 2.3GHz multicore processor. The results are averaged over 20 random runs. The solution of CQP is obtained via SDP relaxation according to \cite{las01,lasspa}, using the software YALMIP, with MOSEK solver, see \cite{yalmip,mosek}. Concerning the SDP relaxation, we observe the relaxation order equal to 2 is generally sufficient to achieve the minimum and the minimizer of the problem.

In Fig. \ref{fig1:errors}, we observe a significant improvement obtained with the proposed CQP with respect to the $\ell_1$ minimization. In particular, when the quantization is sufficiently fine, i.e., more than 1000 quantization levels, the CQP recovery is always exact, while the $\ell_1$ approach slightly improves with the decreasing of the perturbation.  CQP is more accurate, at the price of a slight increase of run time: 9.5 seconds versus 9.2 seconds on average. We also notice that the false positive rate is small for both methods, which is expected by methods that encourage sparsity in the solution. Thus, it is more significant to observe the false negative rate, which shows whether the significant parameters are selected. In Fig. \ref{fig1:errors}, we see that the $\ell_1$ approach never achieves a false negative rate below $0.375$, which means that the erasure of significant parameters is quite frequent.

From Fig. \ref{fig2:errors}, we see that when $\xtrue$ is not binary, the estimate provided by  CQP is still favorable with respect to the $\ell_1$ approach, in terms of relative square $\ell_2$ error and  false negative rate. In particular, we notice that the $\ell_1$ approach does not reach a  false negative rate below $\frac{1}{2}$, which is particularly critical. CQP instead achieves values below $\frac{1}{4}$. On the other hand, the false positive rate is a bit unfavorable for CQP, but generally below $\frac{1}{10}$. We also remark that the false negative rate is more relevant for those  applications where erasing significant entries is more damaging than preserving a few more parameters. In this second experiment, the run time of CQP is slightly increased, passing from 9.5 seconds to 10.5 seconds, which is however acceptable.
\section{Conclusions}\label{sec:con}
In this paper, we address the problem of sparse linear regression from few measurements, when both the regression matrix and the measurements are known in low-precision, or quantized. In the literature, the problem is tackled via $\ell_1$ minimization to promote sparsity on the solution; however, this approach is not reliable when different sparse solutions are present in the feasible set. In this paper, we propose a concave quadratic programming approach which, under suitable conditions, is efficient in the support recovery as well as in the estimation of the non-zero parameters. Sufficient conditions that guarantee the exact recovery of the support are provided. Numerical results show  enhancement with respect to the $\ell_1$ approach, at the price of a slightly increased computational complexity. Future work will envisage a deeper analysis of the performance and the refinement of the sufficient conditions. 

\bibliographystyle{plain}
\bibliography{refs}
\end{document}